\documentclass[12pt,letterpaper,reqno]{amsart}

\usepackage{amssymb}
\usepackage{amsmath,amscd}
\usepackage{color}
\usepackage{epsf}
\usepackage{graphicx}

\theoremstyle{plain}
\newtheorem{thm}{Theorem}[section]
\newtheorem{lem}[thm]{Lemma}

\newtheorem{cor}[thm]{Corollary}

\newtheoremstyle{definition}{7pt plus6.3pt minus6.3pt}{7pt plus3pt minus3pt}%
{\rm}{}{\bf}{}{0.75em}{\thmname{#1}\thmnumber{
#2}\thmnote{\sl\stdspace#3}} \theoremstyle{definition}

\newcommand{\Z}{\mathbb Z}

\begin{document}

\title{Surgery on links with unknotted components and three-manifolds}

\author{Yu Guo$^\dagger$}
\address{$^\dagger$Department of Mathematics, Nanjing University, Nanjing, 210093, P.R.China}
\email{finier@yahoo.com.cn}
\author{Li Yu*}
\address{*Department of Mathematics and IMS, Nanjing University, Nanjing, 210093, P.R.China}
\email{yuli@nju.edu.cn}

\subjclass[2000]{Primary 57M25, 57M27; Secondary 57M05.}

\begin{abstract}
 It is shown that any closed three-manifold $M$ obtained by
integral surgery on a knot in the three-sphere can always be
constructed from integral surgeries on a $3$-component link
$\mathcal{L}$ with each component being an unknot in the
three-sphere. It is also interesting to notice that infinitely many
different integral surgeries on the same link $\mathcal{L}$ could
give the same three-manifold $M$.
\end{abstract}
\maketitle

\section{introduction}

 It is well known that every closed, orientable,
 connected $3$-manifold $M$ can be obtained
by integral surgery on a link in $S^3$. Moreover, one may always
find a surgery presentation for $M$ in which each component of the
surgery link is an unknot (see \cite{L}). For convenience, we use
the word \textit{simple $n$-link} to denote an $n$-component link
with all its components being unknots in $S^3$. Then the minimal
number $\nu(M)$ of the components in all integral simple $n$-link
surgery presentations for $M$ is a topological invariant of $M$,
that is:
\[
  \nu(M) := \text{min}\{ n\, |\, L\: \text{is a simple $n$-link in $S^3$ and we can get\ } M \
  \text{by   }
\]

$\text{\qquad\qquad doing an integral surgery on\ } L \}$\\

For example: $\nu(S^3) =0$ and $\nu(L(p,1))
 = 1$ where $L(p,1)$ is a lens space ($p\geq 2$). However, it is not
 easy to compute $\nu(M)$ in general.
  In particular, let $S^3_K(m)$ denote the $3$-manifold got from integral surgery on
 a knot $K\subset S^3$ with surgery index $m$. Then it is easy to see that
 $\nu(S^3_K(m))\leq u(K)+1$, where $u(K)$ is the unknotting number of
 $K$. But in fact, We can prove the following:

\begin{thm} \label{main}
  For any knot $K\subset S^3$ and any integer $m$, $\nu(S^3_K(m))\leq
  3$, i.e. we can always construct $S^3_K(m)$ by doing an integral
  surgery on a simple $3$-link in $S^3$.
\end{thm}

 \textbf{Remark:} In \cite{AD}, D.Auckly defined a topological invariant called
  \textit{surgery number} of a closed $3$-manifolds.
 By his definition, the surgery number of $S^3_K(m)$ is $1$ for any knot $K$. The
 $\nu(M)$ defined above can be considered as
 another type of surgery number which is more subtle than Auckly's in the
 sense that $\nu(S^3_K(m))$ could be different for different knot
 $K$.\\

 The geometric and topological
 properties of $S^3_K(m)$ have been studied intensively, which
 reveals much topological information of the knot $K$ itself.
 Theorem~\eqref{main} ought to be useful to understand
 the geometry and topology of $S^3_K(m)$ and hence $K$ in the future.

\ \\
\section{Turn knot into simple $2$-links}

 In this section, we will introduce some special operations on a knot diagram
  called \textit{skein-move}. We will see that the skein-move along with
  plane isotopies and the Reidemeister moves can turn any knot diagram on a plane into the diagram of a
  simple $2$-link.\\

  First of all , for any knot $K\subset S^3$, we can use plane isotopy and
  Reidemeister moves to turn any diagram of $K$ into the form that
  all crossings in the diagram are on a short arc of $K$.\\

\begin{figure}
      \includegraphics[width=1.3cm]{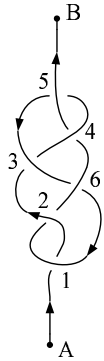}\\
       \caption{Label the crossings of a tangle}\label{p:Label}
   \end{figure}

  The idea is: starting from any diagram $D$ of $K$, we consider $D$ as
  the closure of a $1$-tangle $T$. Then we label the crossings on $T$
  according to their first appearance when we travel
  from the bottom end $A$ of $T$ to the top end $B$ of $T$,
  see figure~\eqref{p:Label} for example. Notice that we
  will meet each crossing of $T$ twice in the process, but when we meet a crossing
  for the second time, we will not relabel it or count it.

\begin{figure}
      \includegraphics[width=10cm]{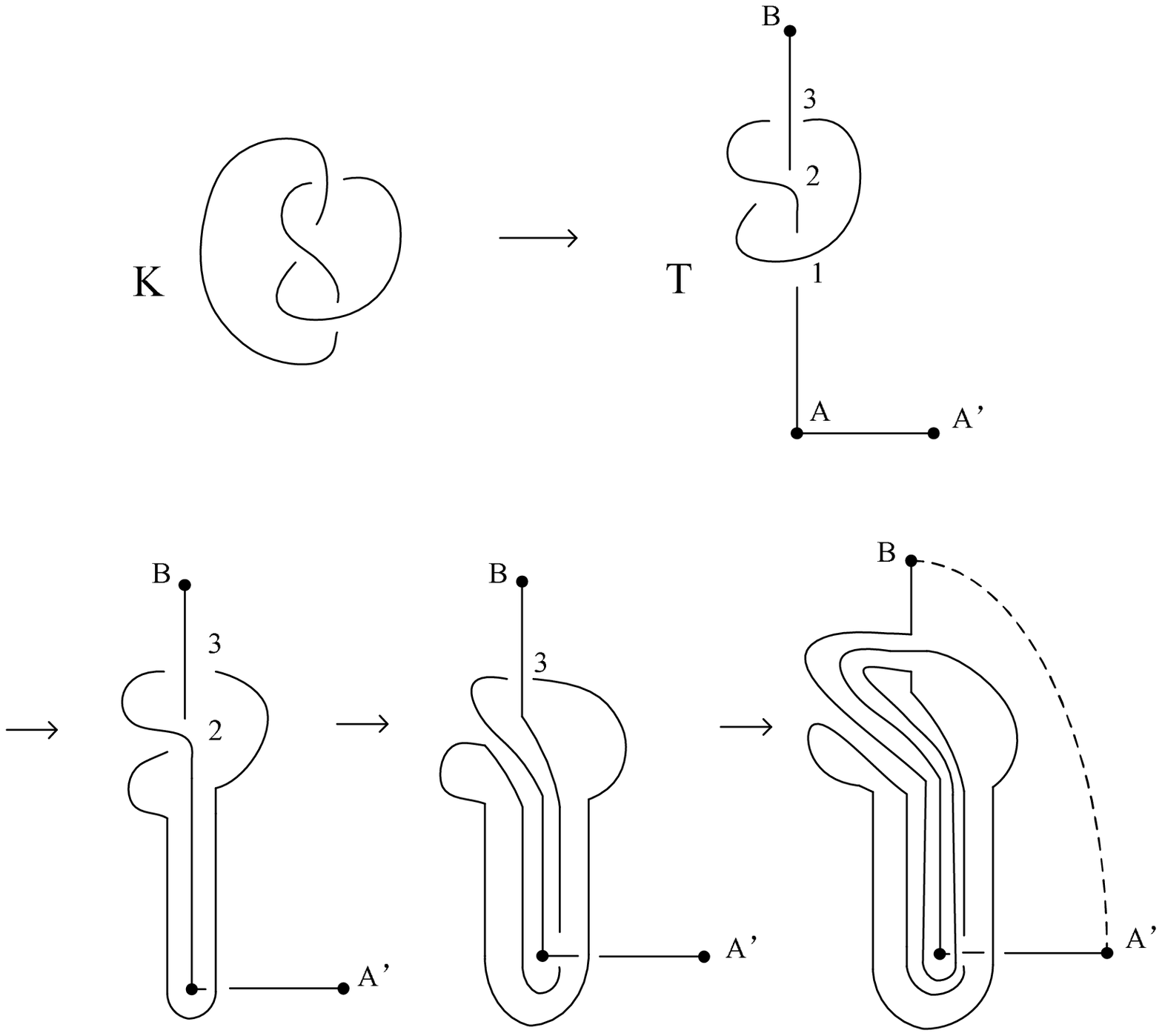}\\
       \caption{tangle}\label{p:tangle}
   \end{figure}

   Next, extend the tangle horizontally via a line
   from $A$ to $A'$.  See the figure~\eqref{p:tangle} for example.

  Denote the crossings by $z_1, \ldots, z_n$ according to their labels.
  Then, start from the crossing $z_1$, we can extend a small segment of
  the strand (overstrand or understrand) at $z_1$ down along the arc of $T$
  that connects $z_1$ and the bottom end $A$, until it meet the line segment
  $\overline{AA'}$. To be more precise,
  when we travel along the tangle starting from $A$ and meet the crossing $z_1$
  at the first time, if we are standing on the understrand of $z_1$, we
  extend the overstrand of $z_1$ down via the process described. Otherwise,
  we extend the understrand of $z_1$ down (See the figure~\ref{p:tangle}).
  Obviously, this will reduce the number of crossings of the tangle above
  the line segment $\overline{AA'}$ by $1$.
  Next, we do the same extension process to strands
  at $z_2, \ldots, z_n$ one by one according to their labeled order.
  When we finish this, all the crossings of the tangle
  will be moved to the segment $\overline{AA'}$.
  Then connect $B$ and $A'$ via a simple arc far away from $T$, we get a
  diagram of $K$ in the required form. This form of knot diagram is called
  \textit{well-posed}.\\

  \textbf{Remark:}
  The Dowker notation (see \cite{A}) of
  a well-posed knot diagram with $m$-crossings
  has the property that: in the two numbers associated to each crossing,
  one is $\leq m$, the other is $\geq m$.\\

\begin{figure}
\centering
\begin{minipage}[c]{0.5\textwidth}
\centering
\includegraphics[width=0.9\textwidth]{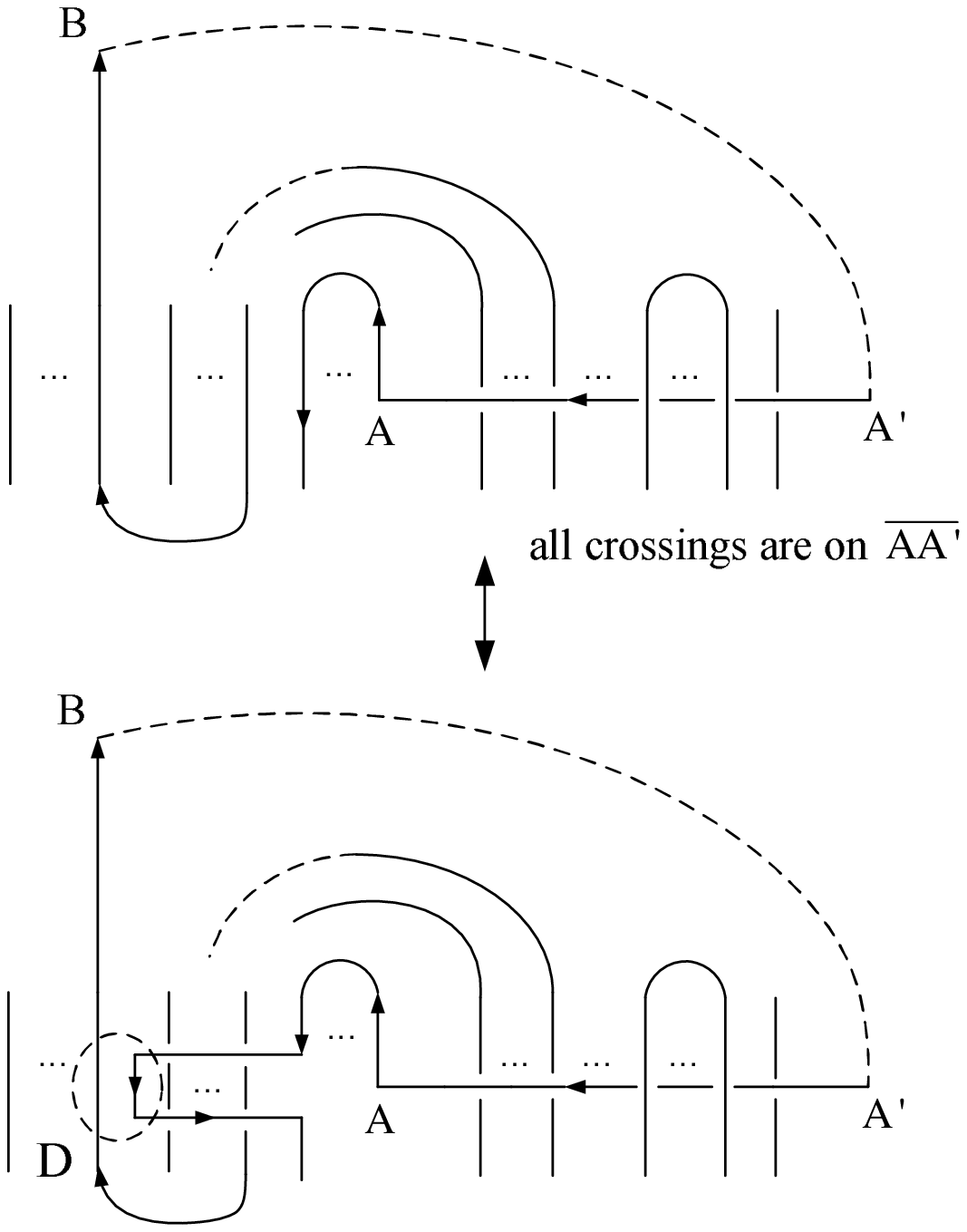}
\end{minipage}%
\begin{minipage}[c]{0.5\textwidth}
\centering
\includegraphics[width=0.9\textwidth]{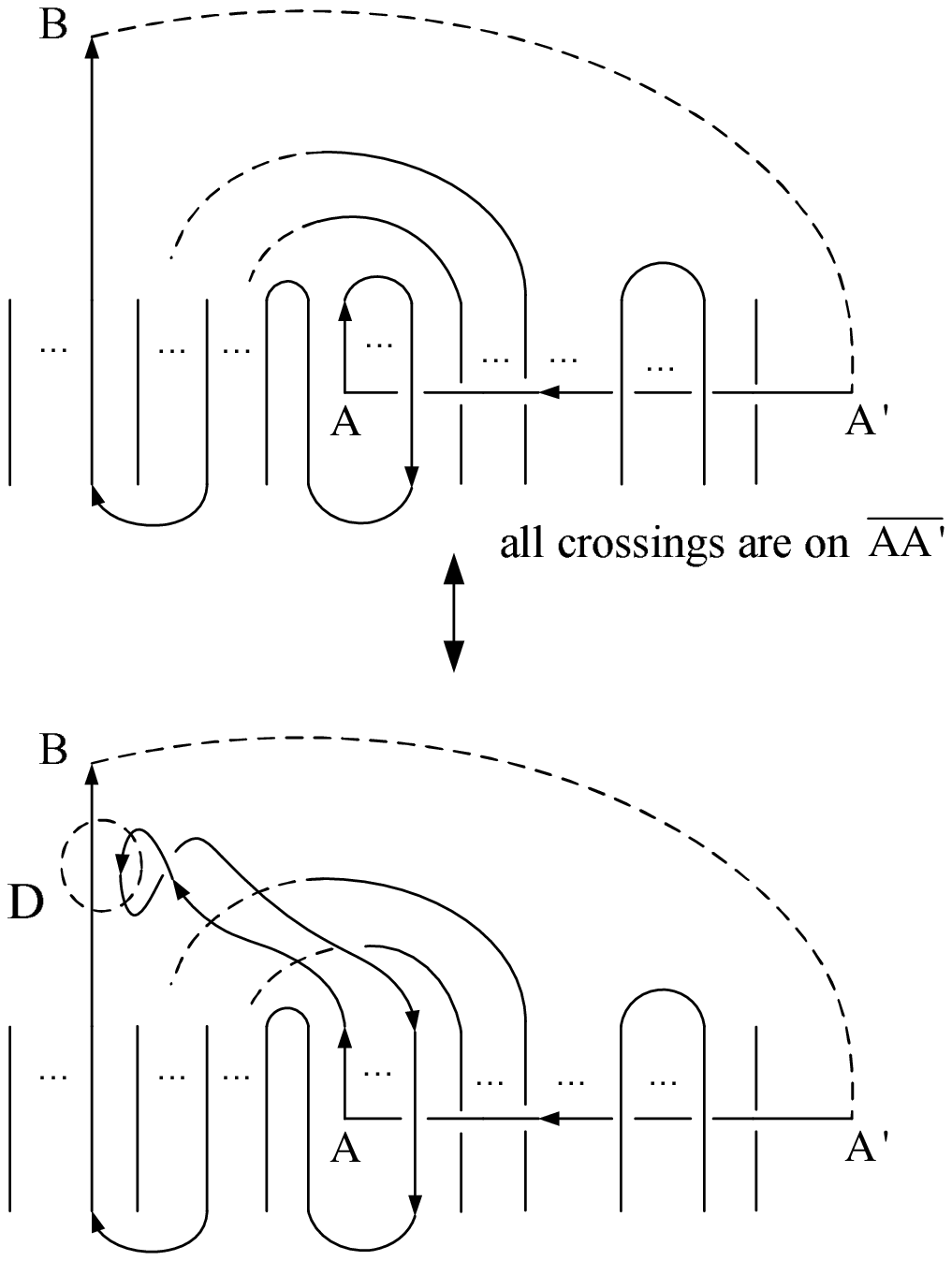}
\end{minipage}\\[20pt]
\begin{minipage}[c]{0.5\textwidth}
\centering
\includegraphics[width=0.8\textwidth]{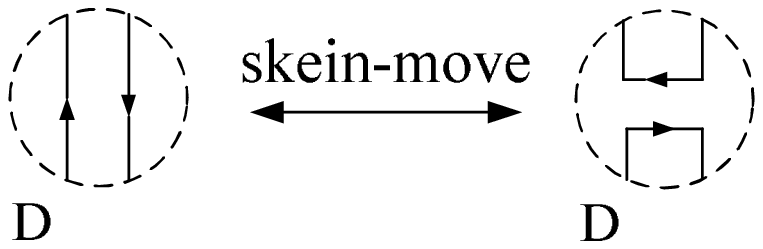}
\end{minipage}%
\caption{ Use skein-move to turn a knot diagram into a simple
$2$-link \label{p:gd}}
\end{figure}

\begin{figure}
  \includegraphics[width=0.7\textwidth]{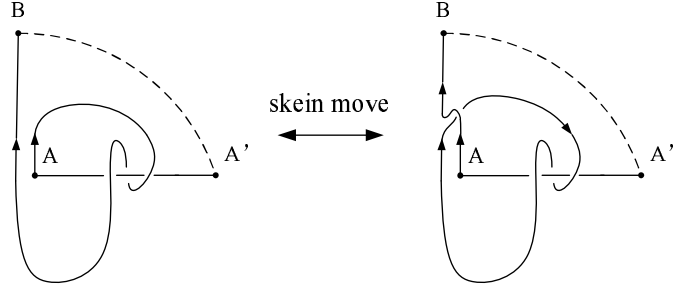}\\
    \caption{Change trefoil knot to a simple $2$-link via a skein move}\label{p:k2l}
\end{figure}

  Next, we orient the knot from $A'$ to $A$. The general picture of a well-posed
  knot diagram is like figure~\ref{p:gd}. Notice that, we can
  always use the skein move defined in figure~\ref{p:gd} to turn
  a well-posed knot diagram into a two-component link $\mathcal{L}$.
  And it is easy to see
  that each component in $\mathcal{L}$ is a diagram of the unknot, i.e.
  the link $\mathcal{L}$ is a simple $2$-link (see figure~\ref{p:k2l} for an
  example).\\

  Conversely, given a diagram of simple $2$-link $\mathcal{L}$, we can use
 Reidemeister moves and the skein move to turn it into a knot
 diagram.\\

\textbf{Remark:} The well-posed diagram for a knot $K$ is not
unique, nor is the corresponding simple $2$-link.\\

\section{$3$-manifolds from integral surgery on a knot}

  Suppose $K$ is a knot in $S^3$, let $N(K)\subset S^3$ be a small tubular neighborhood of
  $K$ and $E(K):=S^3 - N(K)$.  Up to isotopy, $\partial E(K)$ has a \textit{canonical longitude} $l$
  which is homologically trivial
  in $S^3 - K$. And let $m$ be a \textit{meridian} of $\partial E(K)$
  which bounds a disk in $N(K)$. Then doing
   \textit{$(p,q)$-surgery } on $K$ is first removing $N(K)$ from $S^3$ and then
   glue back a standard solid torus $S^1\times D^2$ via a homeomorphism of
   $h: \partial D^2 \times S^1 \rightarrow \partial E(K)$ where $h$ maps the
   $\partial D^2\times {0}$ to a curve on $\partial E(K)$ which is isotopic to
   $p\cdot m + q\cdot l$ on $\partial E(K)$. The $3$-manifold we get is denoted by
   $S^3_K(p,q)$. A $(p,q)$-surgery is called
   \textit{integral} if $q=\pm 1$.  Moreover, $S^3_K(p,q)$ is always an orientable
   $3$-manifold.\\

\textbf{Remark:} We do not need to orient the knot $K$ in the
 surgery since the topological type of $S^3_K(p,q)$ depends only on the knot $K$.\\

 Moreover, we can similarly define surgery on any link $\mathcal{L}\subset
 S^3$. The surgery is called \textit{integral} if the surgery on each
 component of $\mathcal{L}$ is integral.

 \begin{thm}[Lickorish\cite{L62} and Wallace\cite{Wa60}]
   Every closed orientable $3$-manifold can be obtained from $S^3$
  by an integral surgery on a link in $S^3$. Moreover, each component
  of the link can be required to be an unknot in $S^3$.
 \end{thm}

 Integral surgery on a link $\mathcal{L} = L_1\cup \cdots \cup L_m$ decides an integer
 $n_i$ for each component $L_i$ in $\mathcal{L}$, which is called a \textit{framing} of
 $\mathcal{L}$. A link $\mathcal{L}$ with a fixed framing will be called \textit{framed
 link}. So we can also say that any closed orientable $3$-manifolds can be
 got from a surgery
 on a framed link in $S^3$.  \\

 Surgery on different framed links may give the same $3$-manifold.
 Following are two elementary operations on a framed link $\mathcal{L}$ called
 \textit{Kirby moves} (see \cite{Kir78}) which do not change the corresponding
 $3$-manifold.\\

 \begin{description}
   \item[\textbf{ K$1$ Move}] Add or delete an unknotted circle with
    framing $\pm 1$ which belongs to a $3$-ball that does not
    intersect the other components on $\mathcal{L}$.\\

   \item[\textbf{ K$2$ Move}] Slide one component $L_1$ onto another component
   $L_2$. Namely, let $L_2^*$ be a longitude of the tubular
   neighborhood of $L_2$ whose linking number with $L_2$ is the
   framing index $n_2$ of $L_2$. Now replace $L_1$ by
   $L'_1 = L_1 \#_b L_2^*$ where $b$ is any band connecting $L_1$ to
   $L_2^*$ and disjoint from the other components of $\mathcal{L}$.
   The framing of $L'_1$ is $n_1 + n_2 + 2lk(L_1,L_2)$
   where $lk(L_1,L_2)$ is the linking
   number of $L_1$ and $L_2$ in $S^3$ with respect to some orientations of them.
   The rest of the framed link $\mathcal{L}$ remains unchanged.
   To compute $lk(L_1,L_2)$, we orient $L_1$ and $L_2$ in such a way
   that together they define an orientation on $L'_1$. So
   different orientations of $L_1$ and $L_2$ may end up with different
   framed links (see \cite{S}).\\
 \end{description}

\begin{figure}
  \includegraphics[width=0.7\textwidth]{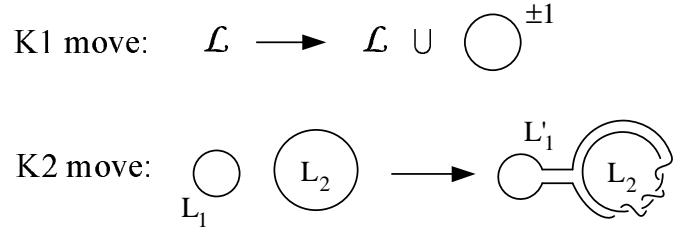}\\
  \caption{Kirby moves}\label{p:kirby}
\end{figure}

  Moreover, it is shown in \cite{Kir78} that any two framed links which give the
  same $3$-manifolds can always be transformed into each other via a finite
  number of Kirby moves. We can use this to show the following
  lemma.

\begin{lem}[proposition $3.3$ \cite{S}] \label{cancel}
If in a framed link $L$ a component $L_0$ is an unknot with framing
zero which links only one other component $L_1$ geometrically once,
then $L_0\cup L_1$ may be moved away from the link $L$ without
changing the resulting $3$-manifold and framings of other
components, and cancelled (See the following figure
\ref{p:cancel})\\
\end{lem}

\begin{figure}
  \includegraphics[width=0.7\textwidth]{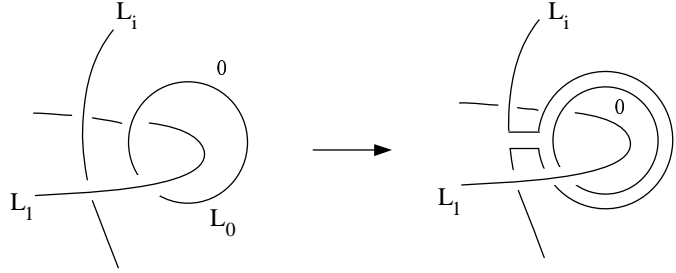}\\
  \caption{move away $L_0$ and $L_1$ then cancel them both}\label{p:cancel}
\end{figure}

 \textbf{Proof of theorem~\ref{main}: }
   By the discussion of previous section, we can turn a knot diagram
   $K$ into a simple $2$-link diagram via plane isotopy, Reidemeister
   moves and the skein move, and vice versa.
   Now suppose $K$ is framed by $m$, it corresponds to
   a closed $3$-manifolds $S^3_K(m)$ via integral surgery.
   Of course, doing skein move to $K$ will not preserve
   the corresponding $3$-manifold in any sense.
    However, the figure~\ref{p:trans} below shows that
   doing an integral surgery on a well-posed diagram of $K$ is equivalent to
   doing integral surgeries on its corresponding simple $2$-link under skein-move
   and an additional $0$-framed unknot. In the figure~\ref{p:trans}, $T$ inside the
   rectangular box represents a tangle of two unknotted strings. The proof is
   just a direct application of Kirby moves and lemma~\ref{cancel}.
    Therefore, integral surgery on a framed
    simple $3$-link can give $S^3_K(m)$ for any knot $K$ and any integer $m$,
   i.e. $\nu(S^3_K(m))\leq 3$.  $\square$\\

   \begin{figure}
  \begin{equation*}
   \vcenter{
            \hbox{
                  \mbox{$\includegraphics[width=\textwidth]{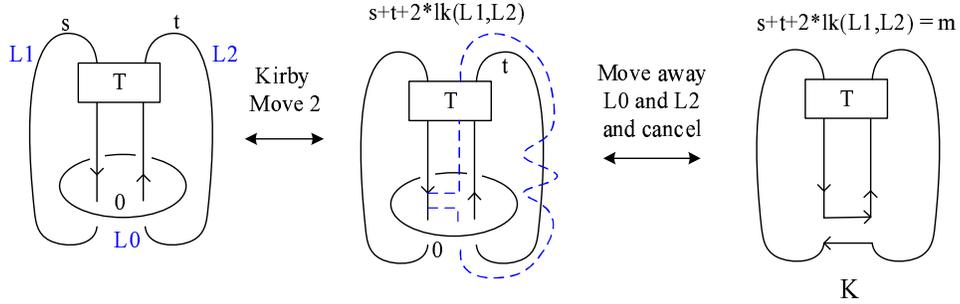}$}
                 }
           }
  \end{equation*}
  \caption{Equivalence of a framed simple $3$-link and a framed knot
      \label{p:trans}}
 \end{figure}

   From the proof of theorem~\ref{main}, we can see the following:
   \begin{enumerate}
    \item  The diagrams for $L_1, L_2$ have no self crossings and
     the geometric intersection number of
   $L_1$ or $L_2$ with a $2$-disk bounded by $L_0$ is $1$.\\

   \item  We can fix the framing on one of the $L_1, L_2$ to be $1$ (or $-1$) in
   the simple $3$-link.\\

   \item There are infinite different framings on a fixed
      simple $3$-link that can give the same $3$-manifold $S^3_K(m)$ !
    \\
    \item We can require the linking number $lk(L_1,L_2)=0$ in the
    simple $3$-link by doing
    second Kirby moves to $L_1$ and $L_0$ in the
    figure~\ref{p:trans}.\\
  \end{enumerate}

  \textbf{Remark:} Obviously, integral surgeries on simple
  $3$-links will give lots of $3$-manifolds other than $S^3_K(m)$. We can
  change the way how $L_0$ is linked to $L_1, L_2$ and
  the surgery index of $L_0$. So theorem~\ref{main} may
   be useful for us to construct some
 interesting examples like integral homology $3$-spheres other than
  $S^3_K(1)$.\\

  \begin{cor}
    Suppose $M$ is constructed from integral surgery on a
    $n$-component link $\mathcal{L}$ in $S^3$, then $\nu(M)\leq 3n$.
  \end{cor}
  \begin{proof}
   Apply the argument in the proof of theorem~\ref{main} to each
   component of $\mathcal{L}$.
  \end{proof}

   Obviously, if $\nu(M^3) = 1$, $M^3$ must be lens space.
   But it is not clear how to classify closed
   $3$-manifolds $M^3$ with $\nu(M^3)=2$. In particular, we can ask
   the following question.\\

   \textbf{Question 1:} For what knot $K$ and integer $m$, $\nu(S^3_K(m)) \leq 2$ ?
      \\

   There are some obvious candidates for the question.
   For example: if the unknotting number of $K$ is $1$,  $\nu(S^3_K(m)) \leq 2$
   for any $m$. But it is not clear how to give a complete answer to this question.
   In particular, it is interesting to know whether
   $\nu(S^3_K(m)) \leq 2$ for all knot $K$ and $m\in \Z$.\\

\begin{figure}
  \includegraphics[width=0.7\textwidth]{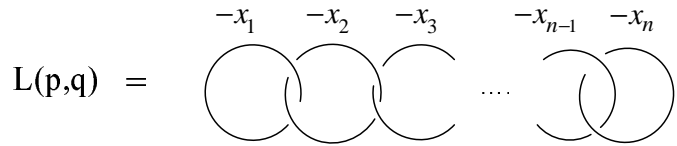}\\
  \caption{surgery diagram for lens spaces}\label{p:lens}
\end{figure}

   Also, it is natural to consider $\nu(S^3_K(p,q))$ for $p\slash q \notin \Z$. For example when $K$
   is the unknot, $S^3_K(p,q)$ is the lens space $L(p,q)$. Suppose the continued fraction decomposition of $p\slash q$
   is $[x_1, \ldots, x_n]$, where
   \[
       [x_1, \ldots, x_n]= x_1 - \frac{1}{\displaystyle x_2
                                - \frac{1}{\displaystyle \cdots x_{n-1}-
                                \frac{1}{\displaystyle x_n}}}
   \]
   then $L(p,q)$ has a surgery presentation as shown in the figure~\ref{p:lens}.
   So $\nu(L(p,q))\leq n$.  Notice that there are examples for $p\slash q =[x_1, \ldots,
   x_n]$ with $n>3$ but $\nu(L(p,q)) \leq 3$. In fact, in \cite{BaRol77}, it is
    shown that
   $L(23,7)$  could be obtained by $-23$-surgery on the $(11,2)$-cable knot about the
   trefoil knot, so $\nu(L(23,7)) \leq 3$ while $23\slash 7 = [4,2,2,3]$.
   More examples of getting lens space via integral surgeries on
   knots in $S^3$ can be found in \cite{BaRol77, FiSt80, BLLi89,Wu90}.\\

  \textbf{Question 2:}
   Does there exist an integer $C$ such that $\nu(L(p,q))\leq C$
   for all $p,q\in \Z$? \\

   \textbf{Remark:} \cite{R2} also gave a way of presenting
   $S^3_K(p,q)$ by an integral surgery on some link.
   But we will not get any universal bounds of $\nu(S^3_K(p,q))$ for all $K$
   and $(p,q)$ via the method in \cite{R2}.\\

  Theorem~\ref{main} provides an interesting way
  to see $S^3_K(m)$ via surgery diagrams. From the proof of theorem~\ref{main},
  we can see that the topological information of $S^3_K(m)$ is completely encoded in
  how $L_1$ and $L_2$ are linked together and the surgery index
  $m$. Notice all the crossings in the diagrams of $L_1\cup L_2$ are between
  $L_1,L_2$. So similar to Dowker notation
  for knots, we can use a sequence of numbers to represent $L_1\cup L_2$.
   This could be interesting in
  its own sense.\\

\end{document}